\documentclass[a4paper]{article}
\usepackage[T2A]{fontenc}
\usepackage[cp1251]{inputenc} 
\usepackage[tbtags]{amsmath}
\usepackage{amsfonts,amssymb,mathrsfs,amscd,comment}
\usepackage{tikz}

\newtheorem{theorem}{Theorem}
\newtheorem{corollary}{Corollary}

\newtheorem{proposition}{Proposition}

\newtheorem{conjecture}{Conjecture}

\newenvironment{remark}
{\noindent{\it Remark\/}.}{\smallskip\par}

\begin{document}

\title{\bf\large\MakeUppercase{Proof of Gal's conjecture for the D series of generalized associahedra}}

\author{Mikhail Gorsky}

\maketitle


In this short note we consider generalized associahedra of type $D.$ We prove that these polytopes are not nestohedra for $n \geq 4,$ but the statement of Gal's conjecture holds for them.

A convex polytope of dimension $n$ is said to be {\it simple} if each of its vertices is contained in precisely $n$ facets. A simple polytope $P^n$ is called a {\it flag} if any set of pairwise intersecting facets $F_{i_1},\ldots, F_{i_k}$ has a  non-empty intersection.

The theory of cluster algebras deals with {\it generalized associahedra}, simple flag polytopes dual to the cluster complexes of finite type algebras (\cite{FZ3}). Each such polytope corresponds in a canonical way to a disjoint union of Dynkin diagrams. Standard associahedra $As^n$ correspond to the diagrams $A_n$; cyclohedra $Cy^n$ correspond to the diagrams $B_n$ and $C_n$. We let $D^n$ denote the series of generalized associahedra corresponding to the diagrams $D_n.$

Let $f_i$ be a number of $i$-dimensional faces of an $n$-dimensional polytope $P.$ As usual, let us introduce
$f(P)(t) = \sum\nolimits_{i=0}^{n} f_i t^i$ and 
$h(P)(t) = f(P)(t-1),$ the $f-$ and $h-$polynomials of the polytope $P^n.$ For simple polytopes the  DehnЦSommerville equations hold (see, e.g., \cite{BP}), these are equivalent to the reciprocity of the $h-$polynomial, that is 
$$h(P)(t)=\sum\limits_{i=0}^{[\frac{n}{2}]} \gamma_i t^i (1+t)^{n-2i}.$$
The polynomial $\gamma(P)(\tau)=\sum\nolimits_{i=0}^{[\frac{n}{2}]} \gamma_i \tau^i$ is called the $\gamma-$polynomial.

\begin{conjecture} \normalfont{(Gal \cite{G})}
\textit{Let $P$ be a flag polytope; then $\gamma_i (P) \geq 0.$}
\end{conjecture}

Consider the set $\mathscr{P}^{\mathrm{cube}}$ of polytopes obtained from the standard cube by succesively shaving off faces of codimension~$2$
Gal's conjecture was proved for flag nestohedra in~\cite{V} using the following two facts: each flag nestohedron lies in $\mathscr{P}^{\mathrm{cube}}$, and the statement of the conjecture holds for all polytopes in $\mathscr{P}^{\mathrm{cube}}.$ The polytopes $As^n$ and~$Cy^n$ are nestohedra (see~\cite{P}), hence the  conjecture is already proven for them. Our main result can be formulated as follows: {\it the polytopes $D^n$ are not nestohedra (for $n \geq 4$), however, they lie in $\mathscr{P}^{\mathrm{cube}}$, therefore Gal's conjecture holds for them.} 

A collection $B$ of non-empty subsets of the set $[n~+~1] = \left\{1,\ldots,n + 1\right\}$ is called a {\it connected building set} on $[n~+~1]$ if: 1) $\left\{i\right\} \in B$ for all $i \in [n~+~1]$ and $[n~+~1] \in B;$ 2) if $S_1, S_2 \in B$ and $S_1 \cap S_2 \neq \emptyset,$ then $S_1 \cup S_2 \in B.$

Let $e_i$ be the vectors of the standard basis of $\mathbb{R}^{n+1}.$ For each $S \in [n~+~1]$ we define the simplex $\Delta_S = \mbox{conv} \left\{ e_i, i \in S \right\},$ where 'conv' denotes the convex hull. The {\it nestohedron} $P_B$ is the Minkowski sum $\sum \Delta_S$ over $S \in B.$ Each nestohedron is a simple polytope (see~\cite{P}).

Let us consider the differential ring $(\mathcal{P}, d)$ of polytopes introduced in \cite{B}, where $d P$ is the formal sum of all the facets of a polytope $P.$ For nestohedra the following formula holds: 
\begin{equation} \label{dpb}
d P_B = \sum\limits_{S \in B / B_{max}} P_{B|_S} \times P_{B / S}
\end{equation}
Here $B_{max}$ is the set of maximal elements of $B$ under the inclusion relation,
$B|_S = \left\{S' \subset S, S' \in B \right\},$ and  $B / S = \left\{ S' = S'' \backslash S, S'' \in B\right\}.$

We need the following two known facts:

(i) The set of facets of $As^n$ corresponds bijectively to the set of diagonals of the regular $(n+3)-$gon, at that two facets intersect if and only if the diagonals corresponding to them do not have common interior points (see \cite{Lee}). 

(ii) The set of facets of $D^n$ corresponds bijectively to the set $Diag \cup D \cup D',$ where
the set $Diag$ consists of the unordered
pairs of centrally symmetric (each to the other) non-diameter diagonals of the regular $2n$-gon; each set of $D$ and $D'$ consists of diameters of the $2n$-gon.
It is convenient to think that
each diameter in $D$ is coloured in one color, and each diameter in $D'$ in another.
Two facets intersect if they correspond to pairs of diagonals which do not have common interior points, or to diameters of the same color, or to diameters of different colours connecting the same antipodal points. (see \cite{FZ3}).

Using these facts, we have the following result.

\begin{proposition} \label{func}
$d D^n = n (\sum\limits_{k = 0}^{n-3} (As^{k} \times D^{n - k - 1}) + 2 As^{n - 1}).$
\end{proposition}


The formula (\ref{dpb}) has an immediate consequence.

\begin{proposition} \label{cond}
Let $P$ be an $n-$dimensional polytope that does not decompose into a nontrivial direct product, and let $P$ have at least $2n+3$ facets which do not decompose into a nontrivial direct product. Then $P$ is not a nestohedron.
\end{proposition}


It is known that the polytopes $As^n$ do not decompose into a nontrivial direct product. Hence, by induction, we obtain the following result from Propositions \ref{func} and \ref{cond}.

\begin{theorem} \label{notnest}
Generalized associahedra of type $D_n$ are not nestohedra for $n \geq 4.$
\end{theorem}


\begin{theorem} \label{dpcube}
Generalized associahedra $D^n$ lie in the set $\mathcal{P}^{cube}$.
\end{theorem}

The proof is by induction on $n:$ it can be shown that $D^n$ can be obtained from $D^{n-1} \times I$ by succesively shaving off faces of codimension 2. It turns out that one can first shave off some $(n-1)$ faces of one of the bases of $D^{n-1} \times 0,$ then a face that is the intersection of two facets (one of which is obtained by shaving, while the other is a face of $D^{n-1} \times I$) corresponding to two diameteres of different colours, and then $(n-3)$ more faces of the same base. Since $D^n$ is a flag polytope, its combinatorial structure is determined by pairwise intersections of its facets. In terms of diagonals of the $2n-$gon one can easily check that the procedure described above leads to the required polytope .

\begin{corollary}
Gal's conjecture holds for generalized associahedra $D^n.$
\end{corollary}

\begin{remark}
By a construction in \cite{BV}, one can use a sequence of shavings of the polytope $D^{n-1} \times I$ leading to $D^n$ to obtain an explicit recurrence relation expressing  $\gamma(D^n)$ in terms of $\gamma(D^k)$ and $\gamma(As^k)$, $k < n.$ See the derivation of a similar relation for the cyclohedron in \cite{BV}.
\end{remark}

The author is grateful to V. M. Buchstaber for suggesting the problem and for his attention to this work, and to V. D. Volodin for valuable discussions.


\end{document}